# Determination of the asymptotic behavior of probabilistic characteristics of arithmetic functions and some other questions of probabilistic number theory

Victor Volfson

ABSTRACT One of the questions of distribution of prime numbers is considered in the article. It is shown what error is obtained from the assumption that the asymptotic density of a sequence of primes is a probability. Various forms of an analogue of the law of large numbers for arithmetic functions and, in particular, the Hardy-Ramunajan theorem are obtained. A method is given for finding asymptotics of the probabilistic characteristics of arithmetic functions.

## 1. INTRODUCTION

It is said [1] that the theory of probability studies not only objects of a random nature, but also deterministic rather complex objects that are poorly described by other methods. The application of probability theory to number theory is studied by probabilistic number theory. We will look at two such complex deterministic objects: primes and arithmetic functions in this work.

Indeed, the conjecture of infinity of twin primes has not yet been proved by methods of number theory, although Hardy and Littlewood gave an accurate estimate of their number using probabilistic methods. The picture is similar with other conjectures about prime numbers. The second chapter of this work will be devoted to one of the questions of the distribution of prime numbers.

An arithmetic function is a function defined on the set of natural numbers and taking values on the set of complex numbers. The name arithmetic function is due to the fact that this function expresses some arithmetic property of the natural series. It is interestimg to study the distribution of arithmetic functions in the natural series. These issues are discussed in the third chapter of this work.

Let's start by considering the concept of the density of a subset of the natural series, which underlies the study of both of these objects.

The density of a subset $A$ of a natural series on an interval $[1.n]$ is determined as follows:

$$\frac{\#A}{n}. \qquad (1.1)$$

Any initial segment of a natural series $\{1,...,n\}$ can be naturally transformed into a probability space $(\Omega_n, \mathcal{A}_n, \mathbb{P}_n)$, by taking $\Omega_n = \{1,...,n\}$, $\mathcal{A}_n$ - all subsets $\Omega_n$, $\mathbb{P}_n(A) = \frac{\#A}{n}$.

Then an arbitrary (real) function $f$ of a natural argument $m$ (or, more precisely, its restriction to $\Omega_n$) can be considered as a random variable $\xi_n$ on this probability space: $\xi_n(m) = f(m)$, $1 \leqslant m \leqslant n$.

Let us compare (1.1) with the definition of the probability measure of a given probability space on the initial segment of a natural series $[1,n]$ and make sure that the density of a subset of a natural series is precisely this finite probability measure:

$$P_n(A) = \frac{\#A}{n}. \qquad (1.2)$$

Now we find the limit of the given probability when the value $n \to \infty$, if it exists (asymptotic density):

$$\lim_{n \to \infty} P_n(A) = \lim_{n \to \infty} \frac{\#A}{n}. \qquad (1.3)$$

It is easy to establish that the given asymptotic density is not a probability, since it does not possess the property of countable additivity [2].

1. ABOUT ONE QUESTION OF DISTRIBUTION OF PRIMES

Let us determine what error we get if we assume that the asymptotic density of a sequence of primes is a probability.

We write the asymptotic law of prime numbers in the form:

$$\pi(n) = \frac{n}{\ln(n)}(1+o(1)).  \qquad (2.1)$$

Having in mind (1.2), (2.1), the density of the subset of primes or the probability of choosing a prime number at random from the interval $[1,n]$ is:

$$P_n = \frac{\pi(n)}{n} = \frac{1}{\ln(n)}(1+o(1)). \qquad (2.2)$$

On the other hand, the probability that a natural number $n$ is divisible by a prime number $p$ is:

$$\frac{1}{n}\left[\frac{n}{p}\right]. \qquad (2.3)$$

Therefore, the probabilities limit (2.3) when the value $n \to \infty$ (asymptotic density), (which based on the above is not a probability) is equal to:

$$\lim_{n \to \infty} \frac{1}{n}\left[\frac{n}{p}\right] = \frac{1}{p}. \qquad (2.4)$$

Let's assume that asymptotic density (2.4) is the probability that a natural number $n$ is divisible by a prime number $p$. Then the probability that a natural number $n$ is not divisible by a prime number $p$ is equal to $1 - \frac{1}{p}$.

Having in mind that under these assumptions, the probability that a natural number $n$ is divided entirely by prime numbers $p_1, p_2, ..., p_l$ is equal to $\frac{1}{p_1 p_2 \cdots p_l}$, that we can consider these events are independent.

Therefore, given that opposite events are also independent, the probability that a natural number $n$ is not divisible by prime numbers $p_1, p_2, ..., p_l$ is equal to:

$$\prod_{i \leq l}(1 - \frac{1}{p_i}).$$

Using the sieve of Eratosthenes, taking into account the above assumptions, the probability that a randomly chosen natural number from an interval $[1,n]$ is prime (in the sense of the density of primes on a given interval) is:

$$\prod_{p \leq \sqrt{n}} (1 - \frac{1}{p}). \qquad (2.5)$$

Based on the Mertens theorem, expression (2.5) is equal to:

$$\prod_{p \leq \sqrt{n}} (1 - \frac{1}{p}) = \frac{2e^{-\gamma}}{\ln(n)} (1 + o(1)). \qquad (2.6)$$

where $\gamma$ is Euler's constant.

Let us compare asymptotics (2.2) and (2.6). The leading terms of these asymptotics differ by a constant value $2e^{-\gamma}$.

2. **VARIOS FORMS ON AN ANALOGUE OF THE LAW OF LARGE NUMBERS AND A METHOD FOR DETERMINATING ASYMPTOTICS OF THE PROBALISTIC CHARACTERISTICS OF ARITHMETIC FUNCTIONS**

Recall that any initial segment of a natural series $\{1,...,n\}$ can be naturally transformed into a probability space $(\Omega_n, \mathcal{A}_n, \mathbb{P}_n)$ by taking $\Omega_n = \{1,...,n\}$, $\mathcal{A}_n$ - all subsets $\Omega_n$, $\mathbb{P}_n(A) = \frac{\#A}{n}$. Then an arbitrary (real) function $f$ of a natural argument $m$ (or, more precisely, its restriction to $\Omega_n$) can be considered as a random variable $\xi_n$ on this probability space: $\xi_n(m) = f(m)$, $1 \leqslant m \leqslant n$.

Therefore, we can write down the Chebyshev inequality for an arithmetic function $f(m), m = 1,...,n$ on a given probability space:

$$P_n(|f(m) - A_n| \leqslant b\sigma_n) \geq 1 - 1/b^2. \qquad (3.1)$$

where $b \geq 1$, and $A_n, \sigma_n$, accordingly, the average value and standard deviation $f(m), m = 1,...,n$.

We put $b = b(n)$ in (3.1), where $b(n)$ is an unboundedly increasing function when the value $n \to \infty$:

$$P_n(|f(m) - A_n| \leq b(n)\sigma_n) \geq 1 - 1/b^2(n). \tag{3.2}$$

The limit $P_n$ in (3.2) with the value $n \to \infty$ exists:

$$P_n(|f(m) - A_n| \leq b(n)\sigma_n) \to 1, n \to \infty. \tag{3.3}$$

Expression (3.3) is an analogue of the law of large numbers for arithmetic functions [3].

Using other forms of Chebyshev's inequality, an analogue of the law of large numbers for arithmetic functions can be written as:

$$P_n(|f(m) - A_n| \geq b(n)\sigma_n) \leq 1/b^2(n). \tag{3.4}$$

or

$$P_n(|f(m) - A_n| \geq b(n)\sigma_n) = O(1/b^2(n)). \tag{3.5}$$

or

$$\#\{m : m \leq n, |f(m) - A_n| \geq b(n)\sigma_n) = O(n/b^2(n))\}. \tag{3.6}$$

The limit in (3.4) when value $n \to \infty$ exists and is equal to:

$$P_n(|f(m) - A_n| \geq b(n)\sigma_n) \to 0, n \to \infty. \tag{3.7}$$

Now, based on the analogue of the law of large numbers for arithmetic functions (3.3), (3.5) and (3.7), we can write the Hardy-Ramunajan theorem [4] in the form (3.8), (3.9), (3.10).

It is known that it is true for the arithmetic function of the number of prime divisors of a natural number $\omega(m), m \leq n$: $A_n = \sigma_n^2 = \ln\ln(n) + O(1)$. Let's take $b(n) = (\ln\ln(n))^\epsilon$, where $\epsilon > 0$, and get:

$$P_n(|\omega(m) - \ln\ln(n)| \leq (\ln\ln(n))^{1/2+\epsilon}) \to 1, n \to \infty. \tag{3.8}$$

$$P_n(|\omega(m) - \ln\ln(n)| > (\ln\ln(n))^{1/2+\epsilon}) \to 0, n \to \infty. \tag{3.9}$$

$$P_n(|\omega(m) - \ln\ln(n)| \geq (\ln\ln(n))^{1/2+\epsilon}) = O(1/(\ln\ln(n))^{1/2+2\epsilon}). \tag{3.10}$$

Let me remind you what is performed for a strongly additive arithmetic function $f(m)$ for each prime $p$: $f(p^a) = f(p)$. $\omega(m)$ is a strongly additive arithmetic function.

Turan proved [5] that if an arithmetic function $f(m), m = 1,...,n$ is strongly additive and for all primes $p$ satisfies the condition: $0 \leq f(p) < c$ and $A_n \to \infty$ when the value $n \to \infty$, then the following form of an analog of the law of large numbers holds:

$$P_n(|f(m) - A_n| \leq b(n)\sqrt{A_n}) \to 1, n \to \infty. \qquad (3.11)$$

and the condition $A_n = \sigma_n^2$ is met.

As we said above, each arithmetic function $f(m), m = 1,...,n$ can be associated with a sequence of random variables $\xi_n$ ("living" in different probability spaces). Each such random variable has a distribution function, so it is natural to consider the question about convergence of the sequence of distribution functions corresponding to the arithmetic function to the limiting distribution function.

Erds and Kats proved in 1939 [6] that the limit distribution function for an arithmetic function $\omega(m), m \leq n$ is normal for any $x \in R$:

$$P_n\left(\frac{\omega(m) - \ln\ln n}{\sqrt{\ln\ln n}} \leq x\right) \to \Phi(x), n \to \infty. \qquad (3.12)$$

where $\Phi(x)$ is the standard normal distribution function.

The indicated authors also generalized (3.12) and proved that if $f(m)$ is a strongly additive arithmetic function and for all primes $p$ is satisfied $|f(p)| \leq 1$ and $\sigma_n \to \infty, n \to \infty$, then the following holds:

$$P_n\left(\frac{f(m) - A_n}{\sigma_n} \leq x\right) \to \Phi(x), n \to \infty. \qquad (3.13)$$

Let's define a random variable $f^{(p)}(m) = f(p)$ for each prime number $p(p \leq n)$ in this case, if $p|m$ - $f^{(p)}(m) = f(p)$ and otherwise $f^{(p)}(m) = 0$.

Then $f^{(p)}(m) = f(p)$ with probability $\frac{1}{n}[\frac{n}{p}]$ and $f^{(p)}(m) = 0$ with probability $1 - \frac{1}{n}[\frac{n}{p}]$.

Therefore, the average value $f^{(p)}(m)$ on the interval $[1,n]$ is:

$$E[f^{(p)}, n] = \frac{f(p)}{n}\left[\frac{n}{p}\right].$$

It is performed for strong additive arithmetic function $f(m)$:

$$f(m) = \sum_{p \leq n} f^{(p)}(m).$$

Therefore, the average value $f(m)$ on the interval $[1, n]$ is:

$$A_n = \sum_{p \leq n} \frac{f(p)}{n}\left[\frac{n}{p}\right].$$

The following asymptotic behavior of the mean value of a strongly additive arithmetic function $f(m)$ implies:

$$A_n \to \sum_{p \leq n} \frac{f(p)}{p}, n \to \infty. \qquad (3.14)$$

Let me remind you that if an arithmetic function $f(m), m = 1, ..., n$ is strongly additive and for all primes $p$ it satisfy conditions: $0 \leq f(p) < c$ and $A_n \to \infty$ when a value $n \to \infty$, then $A_n = \sigma_n^2$ (see 3.11).

Therefore, having in mind (3.13) and (3.14), for a strongly additive arithmetic function $f(m), m = 1, ..., n$ satisfying the conditions: $0 \leq f(p) \leq 1$ and $A_n \to \infty$ when a value $n \to \infty$, it is true:

$$P_n(\frac{f(m) - A_n}{\sqrt{A_n}} \leq x) \to \Phi(x), n \to \infty, \qquad (3.15)$$

where $A_n \to \sum_{p \leq n} \frac{f(p)}{p}, n \to \infty$.

However, the condition $A_n \to \infty$ is not satisfied for all strongly additive functions $f(m), m = 1, ..., n$ that satisfy the condition $0 \leq f(p) \leq 1$.

For example, let's look at the strongly additive arithmetic function

$$f(m) = |\ln \frac{\varphi(m)}{m}|.$$

The function $f(p) = |\ln \frac{\varphi(p)}{p}| = |\ln(1 - 1/p)|$ decreases monotonically with increasing $p$. The maximum value of this function is reached when the value $p = 2$ and is equal to $\ln 2 < 1$, hence the condition $0 \leq f(p) \leq 1$ is met.

On the other hand, considering that $\ln(1 - 1/p) = 1/p + 1/2p^2 + ...$ the series

$$\sum_{p=2}^{\infty} \frac{f(p)}{p} = \sum_{p=2}^{\infty} \frac{|\ln(1 - 1/p)|}{p} \leq \frac{1}{p^2} + \frac{1}{p^3} + ... = \frac{1}{p^2(p-1)} \text{ - converges.}$$

Therefore, the condition $A_n \to \infty$ for a strongly additive arithmetic function $f(m) = |\ln \frac{\varphi(m)}{m}|$ is not met.

Now let's consider a method for determining asymptotic of the probabilistic characteristics of arithmetic functions.

It is often arise problems [7] with the definition of the asymptotic of the mean value of arithmetic functions, and even more in the determination of the asymptotic of the moments of higher orders. Therefore, I propose to approach the solution of this problem from the other side.

The prerequisites for this approach are as follows.

1. Let there are two sequences of random variables that have one limiting distribution. Naturally, the asymptotic of the mean value and the moments of higher orders of these sequences coincide in this case.

2. It is known that an arithmetic function can be represented as a sequence of random variables. If this sequence converges to the limiting distribution function, which coincides with the limiting distribution function of another sequence of random variables, then the arithmetic function and the sequence of random variables have the same asymptotic for the mean value and moments of higher orders. In this case, it is possible to construct a sequence of random variables for which the specified characteristics are determined more simply. On the other hand, this sequence of random quantities must converge to the same distribution function as the arithmetic function. Then, using the asymptotics of the characteristics of a given sequence of random variables, one can more simply determine the asymptotics of the characteristics of an arithmetic function.

3. If two arithmetic functions have the same limiting distribution, then they have the same asymptotic behavior of all moments.

Let's start with the strongly additive arithmetic functions.

Let me remind you that by definition a strongly additive arithmetic function is a function for which $f(p^a) = f(p)$. Therefore, for a strongly additive arithmetic function and an arbitrary natural number $m = p_1^{a_1}...p_t^{a_t}$, we have:

$$f(m) = f(p_1^{a_1}...p_t^{a_t}) = f(p_1^{a_1}) + ... + f(p_t^{a_t}) = f(p_1) + ... + f(p_t) = \sum_{p|m} f(p). \qquad (3.16)$$

We now consider examples of strongly additive arithmetic functions that satisfy conditions (3.15).

1. Suppose that for a strongly additive arithmetic function $f_1(m)$ for all primes $p$ there are $f_1(p) = 1$. The asymptotic of the mean value of the strongly additive arithmetic function in this case is:

$$\sum_{p \leq n} \frac{1}{p} = \ln\ln(n) + O(1). \qquad (3.17)$$

Therefore, conditions (3.15) are satisfied. Formula (3.17) is also true for the asymptotic behavior of the variance $f_1(m)$. $f_1(m)$ is the number of prime divisors of a natural number $m - \omega(m)$.

2. Suppose that stronger additive arithmetic function $f_2(m)$ has the values $f_2(q_j) = 0$ for prime numbers $q_1, q_2, ..., q_l$, and for other primes - $f_2(p) = 1$. Then it is executed - $\sum_{j=1}^{l} f_2(q_j) = 0$. Therefore, the asymptotic of the mean value $f_2(m)$ differs from the asymptotic (3.17) by $O(1)$ and is also determined by this formula.

Hence, the conditions (3.15) are satisfied. Formula (3.17) is also true for the asymptotic behavior of the variance $f_2(m)$. $f_2(m)$ is the number of prime divisors of a natural number $m$, except prime divisors $q_1, q_2, ..., q_l$.

3. Let a strongly additive arithmetic function $f_3(m)$ take values $0 < f_3(p_i) < 1$ for primes $p_1, p_2, ..., p_k$, and for other primes - $f_3(p) = 1$. Then it is true $\sum_{i=1}^{k} f_3(p_i) < c_1$, where $c_1$ is the constant. Therefore, the asymptotic of the mean value $f_3(m)$ differs from the asymptotic (3.17) by $O(1)$ and is also determined by this formula.

Hence, the conditions (3.15) are satisfied. Formula (3.17) is also true for the asymptotic behavior of the variance $f_3(m)$. $f_3(m)$ is the number of prime divisors of a natural value $m$, where instead of the value $f_3(p) = 1$ will be a constant number greater than 0 but less than 1 for prime numbers $p_1, p_2, ..., p_k$, in the sum (3.16).

4. Let a strongly additive arithmetic function $f_4(m)$ take values $f_4(p_{2k+1}) = 0$ for odd-numbered primes $p_1, p_3, ... p_{2k+1}$ and $f_4(m)$ take values $f_4(p_{2k}) = 1$ for even-numbered primes $p_2, p_4, ... p_{2k}$. Then the asymptotic of the mean value is:

$$\sum_{p_{2k}<n} 1/p_{2k} = 0,5 \ln \ln(n) + O(1). \qquad (3.18)$$

Therefore, conditions (3.15) are satisfied. Formula (3.18) is also true for the asymptotic behavior of the variance $f_4(m)$. $f_4(m)$ is the number of prime divisors with even numbers in a natural number $m$.

Let's look at an example of a strongly additive arithmetic function that satisfies the conditions (3.13).

5. Let a strongly additive arithmetic function $f_5(m)$ take values $f_5(p_{2k+1}) = -1$ for odd-numbered primes $p_1, p_3, ... p_{2k+1}$ and $f_5(m)$ take values $f_5(p_{2k}) = 0$ for even-numbered primes $p_2, p_4, ... p_{2k}$.

Then the asymptotic of the mean value $f_5(m)$ is:

$$\sum_{p_{2k+1}<n} 1/p_{2k+1} = -0,5 \ln \ln(n) + O(1). \qquad (3.19)$$

The asymptotic of the variance $f_5(m)$ is:

$$0,5 \ln \ln(n) + O(1). \qquad (3.20)$$

A strongly additive arithmetic function $f(m), m = 1,..., n$ has a limit normal distribution when $n \to \infty$ if conditions (3.13) are satisfied, in the particular case - condition (3.15). Therefore, if we construct a random variable $S_n = \sum_{p \leq n} X_p$ ($X_p$ - bounded, independent random variables), for which the asymptotic of the mean value and the standard deviation will coincide with the corresponding values of the strongly additive arithmetic function $f(m), m = 1,...,n$ when $n \to \infty$:

$$E[S,n] = A_n, \sigma[S,n] = \sigma_n,$$

then the random variable $S_n$ will also have the same asymptotic behavior of the moments of higher orders with the corresponding values $f(m), m = 1,...,n$ when $n \to \infty$.

Let us show this using Example 1 for a strongly additive arithmetic function - the number of prime divisors of a natural number - $\omega(m)$.

It is performed for the function $f(p) = \omega(p) = 1$. Therefore, the asymptotic of the mean value $\omega(m)$ on the interval $[1, n]$ when $n \to \infty$:

$$E[\omega, n] \to \sum_{p \leq n} \frac{f(p)}{p} = \sum_{p \leq n} \frac{1}{p} = \ln\ln(n) + O(1).$$

and row $\sum_{p=2}^{\infty} \frac{\omega(p)}{p}$ - diverges.

As mentioned earlier, the asymptotic behavior of the variance $\omega(m)$ when the value $n \to \infty$ is also equal to:

$\ln\ln(n) + O(1)$.

Thus, all the conditions for convergence to a normal distribution with the indicated characteristics for a strongly additive function $\omega(m)$ are satisfied.

We define the central moments of higher orders for the arithmetic function $\omega(m)$.

Let us construct a sequence of random variables, converging to a normal distribution with similar mean and variance values.

We use the old notation. Let's consider a random variable with a distribution:

$$P(X_p = 1) = 1/p, P(X_p = 0) = 1 - 1/p,,$$

where $p$ is a prime number. Then it is executed $E[X_p] = 1/p, D[X_p] = 1/p - 1/p^2$.

Let $S_n = \sum_{p \leq n} X_p$, then it is executed $E[S_n] = \sum_{p \leq n} 1/p = \ln\ln(n) + O(1)$.

Suppose that the random variables $X_p$ are independent, then it is fulfilled $D[S_n] = \sum_{p \leq n} 1/p - 1/p^2 = \ln\ln(n) + O(1)$, since the series $\sum_{p=2}^{\infty} 1/p^2$ - converges.

A random variable $S_n$ is the sum of independent and bounded random variables; therefore, based on the Central Limit Theorem (CLT) a sequence of random variables $S_1, S_2,...$ converges to a normal distribution.

The arithmetic function $\omega(m)$ is strongly additive and satisfies the condition of convergence to the normal distribution, as shown above.

The limiting distributions of a sequence of random variables $S_1, S_2,...$ and an arithmetic function $\omega(m)$ coincide, since the asymptotics of their mean values and variances coincide, and, consequently, the asymptotics of all moments over high orders.

Determining the asymptotics of moments of higher orders for an arithmetic function $\omega(m)$ from the asymptotics of the corresponding characteristics of a random variable $S_n$ is much easier. Let's show it.

First, we determine the central moment $k$-th order of a random variable $X_p$:

$$E[(X_p - 1/p)^k] = (1 - 1/p)^k 1/p + (0 - 1/p)^k (1 - 1/p) = 1/p + O(1/p^2).$$

Now let's determine the central moment $k$ - th order of a random variable $S_n$:

$$\sum_{p \leq n} 1/p + O(\sum_{p \leq n} (1/p^2)) = \ln\ln(n) + O(1)$$

since the series $\sum_{p=2}^{\infty} 1/p^2$ - converges.

Therefore, the arithmetic function has the asymptotics of the mean value and all central moments of higher orders coincide and they are equal:

$$\ln\ln(n) + O(1) . \tag{3.21}$$

This method can be used to determine the asymptotics probabilistic characteristics of other strongly additive arithmetic functions that satisfy the conditions (3.13) or (3.14).

The class of strongly additive functions is expanded in [3], and the class H of arithmetic functions satisfying the following conditions is investigated.

Let $f(m) \in H$ and $q^b$ is a sequence of positive integer prime numbers such that: $\sum_b 1/q^b$ - converges.

Let us define another additive function $f^*(m)$, setting $f(p^b) = f^*(p^a)$ for all $p^a$ other than $q^b$. Let $f^*(q^b)$ - run through any values for numbers $q^b$. Then it was proved that the limiting laws: $P_n(\frac{f(m) - A_n}{D_n} < x)$ and $P_n(\frac{f^*(m) - A_n}{D_n} < x)$ exist and coincide.

In particular, we can take as a strongly additive arithmetic function $f^*(m)$ associated with $f(m)$ the relation $f^*(p^a) = f(p)(a = 1, 2, ...)$ for all prime $p$, when $n \to \infty$.

Let us consider a method for determining the asymptotics of the probabilistic characteristics of arithmetic functions, which, as was proved in [3], belong to the class $H$

The arithmetic functions of the number of divisors of a natural number, taking into account their multiplicity - $\Omega(m)$. This arithmetic function (for all prime values) coincides with a strongly additive function $\omega(m)$, i.e. $\Omega(p) = \omega(p)$. Therefore, the asymptotics of the mean value and variance for these functions coincide and are equal to $\ln\ln(n) + O(1)$.

Since the arithmetic function $\Omega(m), m = 1, 2, ..., n$ tends to a normal distribution when a value $n \to \infty$ with similar characteristics, as a strongly additive function $\omega(m)$, then all other characteristics of the arithmetic function coincide. Therefore, the asymptotics of all central moments of the arithmetic function $\Omega(m), m = 1, 2, ..., n$ are also equal $\ln\ln(n) + O(1)$.

Arithmetic functions of the number of prime divisors located on the sequences $4k + 1, 4k - 1$, respectively $\omega_1(m), \omega_2(m)$: also have a limiting normal distribution with the asymptotics of mean values and variances which are equal to $0,5\ln\ln(n) + O(1)$.

Let's consider a random variable $X_p$ that takes two values: $X_p = 1$ with probability $P(X_p = 1) = 1/2p$ and $X_p = 0$ with probability $P(X_p = 0) = 1 - 1/2p$. Then the mean is $E[x_p] = 1/2p$ and the variance is $D[X_p] = 1/2p - 1/4p^2$.

Let's construct a random variable $S_n = \sum_{p \leq n} X_p$, where $X_p$ are independent random variables. Then the asymptotics of the mean value and variance of the given random variable are respectively equal to: $E[S,n] = 1/2 \sum_{p \leq n} 1/p = 0,5 \ln \ln(n) + O(1)$,

$D[S,n] = 1/2 \sum_{p \leq n} 1/p - 1/4 \sum_{p \leq n} 1/p^2 = 0,5 \ln \ln(n) + O(1)$, since the series $\sum_{p=2}^{\infty} 1/p^2$ - converges.

Thus, the mean value, variance and their asymptotics for a sequence of random variables $S_1, S_2,...S_n$ coincide with the mean value, variance and asymptotics of arithmetic functions $\omega_1(m), \omega_2(m), m = 1, 2,.., n$ when a value $n \to \infty$.

Since the sequence of random variables $S_1, S_2,...S_n$ (based on the CLT) also tends to a normal distribution, as in arithmetic functions $\omega_1(m), \omega_2(m), m = 1, 2,.., n$ when the value $n \to \infty$, then the limiting distribution functions for them coincide, and, therefore, all other characteristics coincide.

Let us determine the asymptotics of the central moments of higher orders for arithmetic functions $\omega_1(m), \omega_2(m), m = 1, 2,.., n$ when the value $n \to \infty$.

First, we define the asymptotics of the central moments of the $k$-th order of the random variable $X_p$:

$$E[(X_p - 1/2p)^k] = (1 - 1/2p)^k \frac{1}{2p} + (-1/2p)^k (1 - 1/2p) = 1/2p + O(1/p^2).$$

Now let us define the asymptotics of the central moments of the $k$-th order of the random variable $S_n$:

$$\sum_{p \leq n} E[(X_p - 1/2p)^k] = \sum_{p \leq n} 1/2p + O(\sum_{p \leq n} 1/p^2) = 0,5 \ln \ln(n) + O(1), \tag{3.22}$$

since the series $\sum_{p=2}^{\infty} 1/p^2$ - converges.

These characteristics have correspondingly also arithmetic functions $\omega_1(m), \omega_2(m), m = 1, 2, .., n$ when the value $n \to \infty$.

Let's consider one more additive arithmetic function $\omega_1(m) - \omega_2(m)$. It is proved that this arithmetic function has a limit when $n \to \infty$ - normal distribution with the mean value equal to 0 and the variance equal to $0,5 \ln \ln(n) + O(1)$.

Let us determine asymptotics of the central moments of higher orders of the arithmetic function $\omega_1(m) - \omega_2(m), m = 1, 2, .., n$ when the value $n \to \infty$.

Let's consider a random variable $X_p$ that takes two values $X_p = 1/\sqrt{2p}$ with probability $P(X_p = 1/\sqrt{2p}) = 1/2$ and $X_p = -1/\sqrt{2p}$ with probability $P(X_p = -1/\sqrt{2p}) = 1/2$.

Then the average value of this random variable will be $E[X_p] = 0$, and the variance value will be $D[X_p] = E[X_p^2] = 1/2(\frac{1}{2p} + \frac{1}{2p}) = \frac{1}{2p}$.

Let's take as a random variable $S_n = \sum_{p \le n} X_p$ where all $X_p$ are independent. Then the average value of the random variable - $E[S, n] = 0$, and the variance of the random variable $D[S, n] = \sum_{p \le n} \frac{1}{2p} = 0,5 \ln \ln(n) + O(1)$, i.e., equal to the average value and variance of the arithmetic function $\omega_1(m) - \omega_2(m), m = 1, 2, .., n$ when the value $n \to \infty$.

Based on the CLT, a sequence of random variables $S_1, S_2, ...$ tends to a normal distribution when the value $n \to \infty$.

The limiting distributions of a sequence of random variables $S_1, S_2, ...$ and an arithmetic function $\omega_1(m) - \omega_2(m), m = 1, 2, .., n$ coincide, and, therefore, asymptotics of all their characteristics coincide.

Let us determine asymptotics of the central moments of higher orders of the arithmetic function $\omega_1(m) - \omega_2(m), m = 1, 2, .., n$ when $n \to \infty$.

First, we define asymptotics of the central moments the $k$-th order of the random variable $X_p$:

$$E[(X_p)^k] = 1/2((\frac{1}{\sqrt{2p}})^k + (-\frac{1}{\sqrt{2p}})^k).$$

The value is $E[(X_p)^k] = 0$ if $k$ is odd, and the value is $E[(X_p)^k] = (1/2p)^{k/2}$ if $k$ is even.

Now let us determine asymptotics of the central moments of the $k$-th order of the random variable $S_n$.

Value is equal $E[(S_n)^k] = 0$ if $k$ is odd and the value is equal $E[(S_n)^k] = \sum_{p \leq n}(1/2p)^{k/2}$ if $k$ is even.

Therefore, the value - $D[S,n] = \sum_{p \leq n}(1/2p) = 0,5\ln\ln(n) + O(1)$ when $k = 2$.

The value is $E[(S_n)^k] = \sum_{p \leq n}(1/2p)^{k/2} = O(1)$, if $k > 2$, since the series $\sum_{p=2}^{\infty}(1/2p)^{k/2}$ – converges if $k > 2$.

The arithmetic function $\omega_1(m) - \omega_2(m), m = 1,2,..,n$ also has these characteristics when the value $n \to \infty$.

4. CONCLUSION AND SUGGESTIONS FOR FURTHER WORK

The next article will continue to study the behavior of some sums.

5. ACKNOWLEDGEMENTS

Thanks to everyone who has contributed to the discussion of this paper. I am grateful to everyone who expressed their suggestions and comments in the course of this work.